\input amstex
\documentstyle{amsppt}
\NoRunningHeads 
\def\x{\bullet}
\magnification \magstep 1 

\def\Pt{\Bbb P}
\def\Rt{\Bbb R}
\def\Ct{\Bbb C}
\def\Dt{\Bbb D}

\def\sd{\roman{D}}
\def\r{\roman{r}}

\def\Rp{\accent'137Z}
\def\lp{\char'40l}
\def\ap{a\kern -.32em\lower -0.1ex\hbox{\char'030}\kern -.4em\phantom{.}}

\topmatter
\title Residue in intersection homology and $L_p$--cohomology\endtitle
\author Andrzej Weber \endauthor
\thanks Partially supported by KBN 2P30101007\endthanks

\affil Institute of Mathematics, University of Warsaw \endaffil
\address ul. Banacha 2, 02--097 Warszawa, Poland \endaddress
\email aweber\@mimuw.edu.pl \endemail
\abstract We consider a residue form for a singular hypersurface $K$ with 
isolated singularities. Suppose there are neighbourhoods of
the singular points with coordinates in which hypersurface is described by 
quasihomogeneous polynomials. We find a condition on the weights under which
the norm of the Leray residue form is square integrable. For $\dim K\geq
2$ all simple singularities
satisfy this condition. Then the residue form determines an
element in intersection homology of $K$.
We also obtain a residue class in the cohomology of $K$.
\endabstract
\endtopmatter
\document

\head 0. Introduction \endhead

Let $M$ be a complex manifold of dimension $n+1$ and let $K$ be a smooth
hypersurface. Let $Tub\,K$ be a tubular neighbourhood of $K$. 
Consider the diagram:
$$\CD H^*(M\setminus K) @>\delta>>H^{*+1}(M,M\setminus K) @=
H^{*+1}(Tub\,K,Tub\,K \setminus K)\\
@. @V [M]\cap VV @A \tau AA\\
@. H^{BM}_{2n+1-*}(K) @<[K]\cap << H^{*-1}(K)\,.\endCD$$
In the diagram $H^{BM}_*$ denotes Borel--Moore
homology, i.e. homology with closed supports. All coefficients are in $\Ct$.
The residue map 
$$res=\tau^{-1}\circ\delta: H^*(M\setminus K)\longrightarrow H^{*-1}(K)$$
is defined to be the composition of the differential with the inverse of the 
Thom isomorphism.

Suppose $K$ is singular. Then there is no Thom isomorphism, but we can define
a residue morphism
$$res: H^*(M\setminus K)\longrightarrow H_{2n+1-*}^{BM}(K)$$
$$res\,\omega = [M]\cap\delta\omega$$
If $K$ was nonsingular, then this definition would be equivalent to the
previous one since $\xi\mapsto [K]\cap\xi$ is Poincar\'e duality 
isomorphism and
the diagram above commutes. In general there is no hope to 
lift the residue morphism to  cohomology. For $M=\Ct^{n+1}$ the morphism
$res$ is the Alexander duality isomorphism and $[K]\cap$ may be not onto. For 
the
same reason we can't lift 
the residue morphism to intersection homology of $K$.

Let $\omega$ be 
a holomorphic $n+1$--form with a first order pole on $K$. Then a form
$Res\,\omega\in \Omega^n(K\setminus \Sigma_K)$ can be defined; \cite{Le}, see 
\S 1.
We estimate its norm in \S 4. We assume that $K$ has isolated singularities
locally described by quasihomogeneous polynomials.
Let the weights for a singular point be $a_1,\dots ,a_{n+1}$. We prove that if

$$\kappa = \sum_{i=0}^{n+1}{a_i}>1$$
then the norm of $Res\,\omega$ is square integrable (and even
$L_p$--integrable)
for a special choice of a metric. Applying the isomorphism of 
$L_p$--cohomology and intersection homology (see \S3) we conclude that
$Res\,\omega$ 
determines an element in intersection homology and also in cohomology. The
last one may depend on the choice of coordinates. The question of uniqueness
in sheaf theoretic set up is discussed in \S 5.
For $\kappa\leq 1$ there is a way to define an obstruction (higher residue) to

lift the residue class; see the remark in \S6.
There are some questions one should state:

\item{1)} Is it possible to define residue class in intersection homology for
nonisolated general singularities?

\item{2)} How to define the number $\kappa$ or other numerical obstruction to
lift for an arbitrary, possibly nonisolated singularity?

\item{3)} Does the residue class in cohomology depend on the choice of
coordinates?

I was involved in investigating multidimensional residues by Professor Bogdan
Ziemian (see \cite{Zi}).
 Private talks and his hand written notes motivated me to
deal with this subject. I hope that this paper may be
useful in solving partial differential equations. 
I would also like to thank Professors Zbigniew Marciniak, Piotr Jaworski and
Henryk \Rp o\lp \ap dek for help in preparation of this paper.
\vskip 10pt
Contents

\item {\S 1} Residue form,

\item {\S 2} Topology of a neighbourhood of a singular point,

\item {\S 3} $L_p$--cohomology,

\item {\S 4} Local computation,

\item {\S 5} Uniqueness of the lift,

\item {\S 6} Appendix: the $P_8$ singularity.

\head 1. Residue form\endhead 

We recall the Leray method of defining the residue form \cite{Le}.
Let $\omega$ be a smooth closed $k+1$--form on the complement of the set
$$K=\{(z_1, \dots , z_{n+1})\in \Ct^{n+1}:\, s(z_1,\dots ,z_{n+1})=0\}\,.
$$
Suppose that $\omega$ has a first order pole on $K$; i.e. $s\,\omega$ is a 
global
form on $\Ct^n$. Choose local coordinates in which $s$ is the first
coordinate. This can be done outside the singularities of $K$.
Then $\omega$ can be written locally as 
$$\omega={ds\over s}\wedge\r + \eta\,,$$
where $\r$ and $\eta$ do not contain $ds$ and are smooth on $K$. 
Let $\Sigma_K$ be the singular set of $K$.
The form $$Res\,\omega=\r_{|K}\in\Omega^k(K\setminus\Sigma_K)$$ 
is called the residue form of $\omega$. 
The residue form does not depend on the choice of coordinates and on the 
function
describing
$K$. Thus it is defined globally for a hypersurface in a complex manifold.
Moreover, $\r$ is closed on $K$ and its class in 
$H^{k-1}(K\setminus\Sigma_K)$
does not depend on the representant of the
class $[\omega]\in H^k(M\setminus K)$. For smooth $K$ it represents
residue class defined in the
introduction multiplied by $2\pi i$; \cite{SS}, \cite{Do}. In general the 
residue can be defined to be a 
current
on $M$ supported by $K$, \cite{BG}.
If the singularities are isolated ($n>1$) then by Poincar\'e duality for
$K\setminus\Sigma$  we have $H_n(K)=H^n(K\setminus\Sigma)$ and 
$[Res\, \omega] \in H^n(K\setminus\Sigma)$ also coincides with the homology
residue class (multiplied by $2\pi i$) defined in the introduction.

We are particularly interested in holomorphic forms of degree $(n+1,0)$.
Let $\omega$ be such a form. Locally it can be written as
$$\omega={g\over s}dz_1\wedge\dots\wedge dz_{n+1}$$
with $g$ holomorphic. We have
$$ds=\sum_{i=1}^{n+1}{\partial s\over \partial z_i}dz_i\,.$$ 
If ${\partial s\over \partial z_1}\not=0$ then
$$dz_1=\left(\partial s\over \partial
z_1 \right) ^{-1}\left(ds-\sum_{i=2}^{n+1}{\partial s\over \partial 
z_i}dz_i \right) $$ 
and 
$$\omega={g\over s}\left(\partial s\over \partial
z_1\right)^{-1}\left(ds-\sum_{i=2}^{n+1}{\partial s\over \partial
z_i}dz_i\right)\wedge dz_2\wedge\dots\wedge dz_{n+1}=$$ $$
={ds\over s}\wedge g\left(\partial s\over \partial
z_1\right)^{-1} dz_2\wedge\dots\wedge dz_{n+1}$$
To see how $Res\, \omega = g\left(\partial s\over \partial
z_1\right)^{-1} dz_2\wedge\dots\wedge dz_{n+1}$ behaves in the
neighbourhood of the singularities let us calculate its norm in the metric
coming from the coordinate system:
$$|Res\, \omega|_K=\left|{ds\over |ds|}\wedge Res\, 
\omega\right|=\left|{ds\over |ds|}\wedge g\left(\partial s\over \partial
z_1\right)^{-1} dz_2\wedge\dots\wedge dz_{n+1}\right|=$$
$$=\left|{{\partial s \over \partial z_1}\over |ds|}dz_1\wedge 
g\left(\partial s\over \partial
z_1\right)^{-1} dz_2\wedge\dots\wedge dz_{n+1}\right|
={g\over |ds|}\,.$$
We see that $Res\, \omega$ has a pole in singular points.

Suppose $K$ has isolated singularities. Define $K^\circ$ to be $K$ minus the
sum of small balls centered at the singular points of $K$.
Let $j:(K^\circ,\emptyset)\longrightarrow (K^\circ,\partial K^\circ)$ and
$k:\partial K^\circ\longrightarrow K^\circ$ be the inclusions. Then for
$dim\,K=n>1$ we have \cite{Bo}:
$$IH^{\underline m}_n(K)
= im\left(j_*:H_n(K^\circ)\longrightarrow
H_n(K^\circ,\partial K^\circ)\right)
=$$ $$=im \left(j^*:H^n(K^\circ,\partial K^\circ)\longrightarrow
H^n(K^\circ)\right)
=$$ $$= ker\left(k^*:H^n(K^\circ)\longrightarrow H^n(\partial K^\circ)\right)
$$
The morphism $IH^{\underline m}_n(K)\longrightarrow H_n(K)\simeq H^n(K^\circ)$
 is just the
inclusion. Each class $\alpha\in H_n(K)$ is determined by a smooth form
on the nonsingular part of $K$. The group $IH^{\underline m}_n(K)$ consists of
the classes which can be represented by forms with square integrable norms;
see \S3. 
Our goal in \S4 will be to check whether $|Res\,\omega|$ is square integrable,
but first consider the examples.

\example{Example 1.1} Let $s=xy$ and let $\omega={1\over s}dx\wedge dy$. Then
$ds=y\,dx+x\,dy$. The residue form is $dy\over y$ for $x=0$ and $dx\over x$
for $y=0$. 
We see that 
$$K^\circ=(\Ct\setminus B_\epsilon)\times\{0\}\cup \{0\}\times(\Ct\setminus
B_\epsilon)$$ 
and $Res\, \omega$ does not belong to $ker\,k^*=IH^{\underline m}_n(K)$ since
the form $dy\over y$ (and $dx\over x$) is a generator when restricted to the
small circle.\endexample

Since one may think, that the example is degenerated ($K$ is not normal and
$\dim K=1$) let us consider another example. 

\example{Example 1.2}
Let $s$ be a singularity of the type $P_8$:
$$s(z_1,z_2,z_3)=z_1^3+z_2^3+z_3^3\,.$$
The residue class $Res(\frac1sdz_1\wedge dz_2\wedge dz_3)$ has no lift to
intersection homology; see the Appendix.
\endexample

To show that $[Res\,\omega]\in IH^{\underline m}_n(K)=
ker\left(k^*:H^n(K^\circ)\longrightarrow H^n(\partial
K^\circ)\right)$ one has to integrate the residue form over each
$n$-cycle $X\subset \partial K^\circ=K\cap S_\epsilon$. Since the mapping $f$
restricted to $S_\epsilon$ is a fibration in a neighbourhood of
$f^{-1}(0)\cap S_\epsilon$ one can find a continuous family of
cycles $X_t\subset f^{-1}(t)\cap S_\epsilon$ with $X_0=X$. The
function $t\mapsto \int_{X_t}s\omega /ds$ is holomorphic
(single--valued) function. It can be expanded in a series
$\Sigma_\alpha a_\alpha t^\alpha$.
By \cite{ArII p.261} we can assume $-(\alpha+1)\leq d$ where $d$ is
the distance of the Newton diagram of $s$ (see \cite{ArII p. 140}).
In the case of quasihomogeneous $s$ with weights $a_1,\dots
,a_{n+1}$ (see \S4) we have $d=-({a_1}+\dots + {a_{n+1}})$. Thus 
to
deduce that $\int_X Res\, \omega =\int_{X_0}s\omega /ds =0$ one
should take $s$ with ${a_1}+\dots + {a_{n+1}}>1$.  This is 
exactly
the condition obtained in \S 4 by applying $L_p$--methods. 
If $\alpha =0$ occurs in the Taylor expansion of $\int_{X_t}s\omega /ds$ then
0 is a spectral number of the  singular  point.  Properties  (and 
definition) of
spectra were discussed in papers of Varchenko e.g.  \cite{ArII},  \cite{Va} 
and Steenbrink e.g. \cite{St1}, \cite{St2}. This way we have:

\proclaim{Conclusion 1.3} If $\;0\;$ does not belong to the specta of the 
singular
points of $K$ then each residue class lifts to the intersection homology of
$K$.\endproclaim 

From the topological point of view this statement can be partially explained
by the fact that spectral numbers  multiplied  by  $2\pi  i$  are 
logarithms of the
eigenvalues of the monodromy. Thus if 0 is not in the spectrum, then 1 is not
an eigenvalue of the monodromy and K is a rational manifold; see the next
paragraph. 
To check that 0 is not a spectral number for the most of singularities see the
table \cite{ArII, page 275}. This shows that in general a residue class lies 
in
intersection homology.
The $L_p$--method of lift presented in the rest of the paper is not so general
but in addition we obtain a concrete lift of residue class to the cohomology 
of
$K$.

\head 2. Topology of a neighbourhood of a singular point \endhead

Let us assume that $0\in\Ct^{n+1}$ is an isolated singular point of a 
hypersurface
$K$. Intersect $K$ with a ball of a small radius. Then the set
$L=S_\epsilon\cap K$ is called the link of the singular point. Milnor 
\cite{Mi}
gave the precise description of the topology of $L$. It is $2n-1$ dimensional
manifold with nonzero homology only in dimensions 0, $n-1$, $n$ and $2n-1$.
Let $h_*$ be the monodromy acting on the homology of the Milnor 
fiber and let
$\Delta(t)$ be its characteristic polynomial. 

\proclaim{Theorem 2.1. \cite{Mi}, \cite{Hi}} Let $n>2$. The link of an 
isolated singular 
point of $s:\Ct^{n+1}@>>> \Ct$ is
homeomorphic to a sphere if and only with $\Delta (1)=\pm 1$.
The link is a rational homology sphere if and only if $\Delta (1)\not= 0$,
i.e. 1 is not a eigenvalue of the monodromy.
\endproclaim

Milnor describes in his book a recipe for computing $\Delta (t)$ of
quasihomogeneous polynomials. 
We restrict our attention to the case of simple and unimodal 
parabolic
(simply elliptic) singularities. All these types may be represented by
quasihomogeneous polynomials. Our choice is motivated by the following:

\proclaim{Theorem 2.2. \cite{ArI}} Every singularity is simple (i.e. it is of 
the 
type: 
$A_k$,
$D_k$, $E_6$, $E_7$, $E_8$) or it is adjacent to one of the unimodal parabolic
type (i.e. to $P_8$, $X_9$ or $J_{10}$).\endproclaim

Now we list the families of simple singularities and the corresponding 
charac\-te\-ris\-tic
polynomials. The table contains answers to the following questions:
\item {a)} Is the link homeomorphic to a sphere? (For $n=2$ ---
is it a homology sphere?)

\item {b)} Is it a rational sphere?

\settabs 
\+***************************&**** & **n** & characteristic polynomial*** &
yes* & yes \cr 
\+Singularity type &\hfil $k$ & \hfil $n$ & characteristic polynomial
& \hfil a) & \hfil b)\cr
\vskip 5pt
\+$A_k:\,z_1^{k+1}+\sum_{i=2}^{n+1}z_i^2$ & odd \quad& odd 
& $\pm(t^k-t^{k-1}+\dots\pm 1)$&no&no\cr
\+ & even \quad & odd & &yes&yes\cr
\+ & all \quad & even 
& $t^k+t^{k-1}+\dots+1$&no&yes\cr
 
\+$D_k:\, z_1^2z_2+z_2^{k-1}+\sum_{i=3}^{n+1}z_i^2$ & $\geq 4\quad $& odd 
& $\pm(t-1)(t^{k-1}-(-1)t^k)$&no&no\cr
\+ & $\geq 4$\quad & even 
& $\pm(t+1)(t^{k-1}+1)$&no&yes\cr
 
\+$E_6:\,z_1^3+z_2^4+\sum_{i=3}^{n+1}z_i^2$ && odd 
& $t^6-t^5+t^3-t+1$&yes&yes\cr
\+ && even 
& $t^6+t^5-t^3+t+1$&no&yes\cr
 
\+$E_7:\,z_1^3+z_1z_2^3+\sum_{i=3}^{n+1}z_i^2$ && odd 
& $-(t-1)(t^6+t^3+1)$&no&no\cr
\+ & & even 
& $-(t+1)(t^6-t^3+1)$&no&yes\cr
 
\+$E_8:\, z_1^3+z_2^5+\sum_{i=3}^{n+1}z_i^2$ && odd 
& $t^8-t^7+t^5-t^4+t^3-t+1$&yes&yes\cr
\+ && even 
& $t^8+t^7-t^5-t^4-t^3+t+1$&yes&yes\cr

The unimodal parabolic singularities as follows
\footnote{
The number $a$ is such that: $a^3+27\not= 0$ for $P_8$, $a^2\not= 4$ for $X_9$
and $4a^3+27\not= 0$ for $J_{10}$\hfill \ }:
\settabs 
\+********************************** & **n** & characteristic polynomial & 
yes* & yes \cr 
\+Singularity type
 & \hfil $n$ & characteristic polynomial &
 \hfil a) & \hfil b)\cr
\vskip 5pt
\+$P_8:\, z_1^3+z_2^3+z_3^3+az_1z_2z_3+\sum_{i=4}^{n+1}z_i^2$ & odd 
& $(t^3+1)^2(t^2-t+1)$&no&yes\cr
\+ & even 
& $(t^3-1)^2(t^2+t+1)$&no&no\cr
 
\+$X_9:\, z_1^4+z_2^4+az_1^2z_2^2+\sum_{i=3}^{n+1}z_i^2$ & odd 
& $-(t^4-1)^2(t-1)$&no&no\cr
\+ & even 
& $-(t^4-1)^2(t+1)$&no&no\cr
 
\+$J_{10}:\,z_1^3+z_2^6+az_1^2z_2^2+\sum_{i=3}^{n+1}z_i^2$ & odd 
& $(t^6-1)(t^3+1)(t-1)$&no&no\cr
\+ & odd 
& $(t^6-1)(t^3-1)(t+1)$&no&no\cr
 
We see that the link of a singular point often happens to be a rational
homology sphere. If it is the case then $K=\{s=0\}$ is a rational homology
manifold and the Poincar\'e duality map 
$$PD:H^k(K;\Ct)\longrightarrow H_{2n-k}(K;\Ct)$$
is an isomorphism. Thus each residue class lifts to cohomology. For other 
cases
there is no chance to construct the uniform lift of the residue morphism. We
will study only the residues of meromorphic forms with a first order pole on 
$K$.
 
\head 3. $L_p$--cohomology\endhead

To show that the residue form on the nonsingular part of $K$ determines an
element in intersection homology we apply the isomorphism which
was suggested in \cite{BGM}: 
 
\proclaim{Theorem 3.1. \cite{Ch}, \cite{We1}} Let $X$ be a pseudomanifold 
equipped 
with a Riemannian
metric on the nonsingular part. Assume that this metric is concordant with a
conelike structure of the pseudomanifold. Then $H^*_{(p)}(X_0)$, the
$L_p$--cohomology of the nonsingular part, is isomorphic to the
intersection homology with respect to the perversity which is the largest
perversity strictly smaller then the function $F(i)={i\over p}$.
\endproclaim 

Concordance with the conelike structure
means that each singular point has a neighbourhood which is quasiisometric to
the metric cone over the link, i.e. to $cL_x=L_x\times [0,1]/L_x\times \{0\}$
with the metric $t^2dx^2 + dt^2$. The intersection homology of a
pseudomanifold $K$ with isolated singularities is either $H^{2n-*}(K)$ or 
$H^{BM}_*(K)$
or the image of the Poincar\'e duality map  
$im(PD:H^{2n-*}(K)@>>>H^{BM}_*(K))$.
The case depends on the value of the perversity for $2n$. If the dimension is
one then we should take the normalization of $K$ instead of $K$.
The perversity associated with $p\in [2,2+{2\over n-1})$ is the middle
perversity $\underline m$ 
and $\underline m(2n)=n-1$. Thus we obtain:
 
\proclaim{ Corollary 3.2} If a hypersurface $K$ with isolated singularities is
equipped with a conelike metric then 
$$H^n_{(p)}(K\setminus\Sigma_K)\simeq\left\{
\aligned
H^{BM}_n(K)\quad &\text{for}\quad 1+{1\over 2n-1}\leq p< 2\\
im\,PD\qquad &\text{for}\qquad 2\leq p < 2+{2\over{n-1}}\\
H^n(K)\qquad &\text{for}\qquad 2+{2\over{n-1}}\leq p\,.
\endaligned\right.$$\endproclaim

In \S4 we construct a suitable conelike metric and estimate the norm of a
residue form for every $p>1$. We show that it is $L_p$--integrable for a 
wide class of isolated singularities including all simple singularities. In
this way we obtain a lift of the residue class to intersection homology.
 
\head 4. Local computation \endhead
 
Recall that we say that a polynomial is quasihomogeneous with weights $a_1,
\dots ,$ $ a_{n+1}$ $a_i>0$, if it is a sum of monomials $\Pi z_i^{k_i}$ such
that $\sum_{i=1}^{n+1}{{k_i}{a_i}}=1$. The homogeneous polynomial of 
degree $d$ is quasihomogeneous with weights $a_i=\frac 1d$.
We show the following:
 
\proclaim{Theorem 4.1} Let $s$ be a quasihomogeneous polynomial in $n+1$ 
variables 
with
weights 
$a_1,\dots,a_{n+1}$. Assume that 0 is an isolated critical point of $s$. 
If $$\kappa=\sum_{i=2}^{n+1}{a_i}>1$$ 
then
there exists a conelike metric on $\Ct^{n+1}$ such that the norm of
the residue form $$Res\left({g\,dz_1\wedge\dots\wedge dz_{n+1}}\over
s\right)\in \Omega^{n,0}(\{s=0\}\setminus \{0\})$$
is $L_p$--integrable.
\endproclaim

\demo{ Proof} We choose $m\in \Rt$ and parametrize $\Ct^{n+1}$ by the
homeomorphism: 
$$(u_1,\dots ,u_{n+1})\longmapsto 
(u_1|u_1|^{{m\,{a_1}}-1},\dots ,u_{n+1}|u_{n+1}|^{{m\, 
a_{n+1}}-1})\,.$$ 
The set $\Phi^{-1}(K)$ is conical. We estimate the norm of the
residue form 
$$\r =Res\left({g\,dz_1\wedge\dots\wedge dz_{n+1}}\over
s\right)={\left({g\,dz_2\wedge\dots\wedge dz_{n+1}}\over
{\partial s\over dz_1}\right)}$$
in the metric induced by this parametrization.
The norm $|\Phi^*(dz_i)|_u$ is (real) homogeneous of degree
${m\, a_i}-1$, the denominator ${\partial s\over \partial 
z_1}\left(\Phi(u)\right)$
is homogeneous of degree $m-{m\, a_1}$.
Thus the norm $|\Phi^*\r |_u$ is bounded by a homogeneous function of degree
$$\sum_{i=2}^{n+1}\left({m\, a_i}-1\right)-\left(m-{m\, a_i}\right)
=\sum_{i=1}^{n+1}{m\, a_i}-n+m
=m\left[\sum_{i=2}^{n+1}{a_i}-1\right]-n\,.$$
Then
$\int_{\{|u|=r\}\cap K}\,|\Phi^*\r |^p_udz $ is bounded by a 
homogeneous function of degree 
$$\alpha=p\left\{m\left[\sum_{i=2}^{n+1}{
a_i}-1\right]-n\right\}+2n-1= p\,m(\kappa-1)+(2-p)n-1$$
If $p=2$ then we see that this function is integrable. For $p>2$ 
one must take $m$ large enough so that $\alpha >-1$.
\qed\enddemo

\proclaim{ Corollary 4.2} If $K$ has 
quasihomogeneous singularities with $\kappa >1$ then
the residue form defines an element in $L_p$--cohomology of $K$ for a suitable
metric. \endproclaim 

\demo{Proof} In each singular point we choose $m$ such that
$m(\kappa-1)>(p-2)n$. Then the residue 
form is $L_p$--integrable with respect to the conelike metric constructed in
the proof of the Theorem 4.1. Hence it defines an element in
$L_p$--cohomology. 
\qed\enddemo 

\remark{Remarks} In the proof of the Theorem 4.1. we can use the function
$e^{- {a_i\over |z_i|}}$ as well as $|z_i|^{m\, a_i}$ ($m$ large). As a
result we get the same condition for weights.
Note that this condition is fulfilled if the Hessian of $s$ is of rank at
least 
2 and $n\geq 2$. Then $s$ has  either a term $z_iz_j$ or $z_i^2+z_j^2$ so
${a_i}+{a_j}=1$ and the other summands in $\kappa$ are nonzero. Practically
the theorem shows that we can integrate the residue form over chains which are
regular enough i.e. which enter singular points along the cone lines.
\endremark
Below we list singularities with computed weights and with the $\kappa$
numbers. They coincide with 'the oscillation indicators' from \cite{ArII}.
\settabs
\+***************&*******&************************* & ******* \cr 
\+ &Type & weights & $\kappa$ \cr
\vskip 5pt
\+ &$A_k$ & $\frac 1{k+1},\ \frac 12,\ \frac 12,\dots$ &${n\over 2}+{1\over 
k+1}$\cr
\+ &$D_k$ & ${{k-2}\over {2k-2}},\ \frac 1{k-1},\ \frac 12,\ \frac 12,\dots$ 
&${n\over 2}+{1\over 
2(k-1)}$\cr
\+ &$E_6$ & $\frac 13,\ \frac 14,\ \frac 12,\ \frac 12,\dots$ &${n\over 
2}+{1\over 12}$\cr
\+ &$E_7$ & $\frac 13,\ {2\over 9},\ \frac 12,\ \frac 12,\dots$ &${n\over 
2}+{1\over 18}$\cr
\+ &$E_8$ & $\frac 13,\ \frac 15,\ \frac 12,\ \frac 12,\dots$ &${n\over 
2}+{1\over 30}$\cr
\vskip 5pt
\+ &$P_8$ & $\frac 13,\ \frac 13,\ \frac 13,\ \frac 12,\ \frac 12,\dots$ 
&${n\over 2}$\cr
\+ &$X_9$ & $\frac 14,\ \frac 14,\ \frac 12,\ \frac 12,\dots$ &${n\over 2}$\cr
\+ &$J_{10}$ & $\frac 13,\ \frac 16,\ \frac 12,\ \frac 12,\dots$ &${n\over 2}$

\cr
\vskip 5pt
We see that for all simple singularities we have $\kappa>1$ for $n\geq 2$. For
unimodal parabolic singularities one should take $n\geq 3$. The Example 1.2.
(see Appendix) shows,
that for $P_8$, $n=2$, the residue class has no lift to the intersection
homology. 
The lift of $Res\,\omega$ to cohomology we may cal a
regularization, that is giving a meaning to the symbol $\int_X
Res\,\omega$ where $X$ is a cycle intersecting singularities of
$K$. Certain regularization of residue form was described in \cite{Zi}.

\head 5. Uniqueness of the lift\endhead

Denote by $i$ the inclusion $M\setminus K \hookrightarrow M$.
Let $\Omega^\x_{M\setminus K}$ be the sheaf of complex--valued forms on
$M\setminus K$. 
It is a soft resolution of the constant sheaf $\Ct_{M\setminus K}$. Thus
$Ri_*\Ct_{M\setminus K}=i_*\Omega^\x_{M\setminus K}$. The inclusion $i$
induces a 
distinguished triangle.
$$\matrix
\Ct_M \quad &\longrightarrow & Ri_*\Ct_{M\setminus K}
&=& i_*\Omega^\x_{M\setminus K}\,.\\
_{+1}\nwarrow & & \swarrow & & \\ 
& R\Gamma_K\Ct_M& & & \endmatrix$$
By taking the cohomology we get the long exact sequence of the pair 
$(M,M\setminus K)$.
The stalk of $R\Gamma_K\Ct_M$ is:
$${\Cal H}^j_x(R\Gamma_K\Ct_M)\simeq H^j(B_x,B_x\setminus K)
@<[B_x]\cap<\simeq<
H_{2n+2-j}^{BM}(K\cap B_x)\,,$$
where $B_x$ is a small ball around $x$. Moreover,
the whole sheaf $R\Gamma_K\Ct_M$ is isomorphic (with a shift of degrees)
to the dualizing sheaf: $$ R\Gamma_K\Ct_M[2n+2]\simeq \Dt_K\,.$$
Thus we get the residue morphism (Grothendieck residue)
$$res:\,i_*\Omega_{M\setminus K}^\x[2n+1]= Ri_*\Ct_{M\setminus 
K}[2n+1]@>{+1}>>
R\Gamma_K\Ct_M[2n+2]\simeq \Dt_K$$
which is an isomorphism on the cohomology sheaves for $j\not= -(2n+1)$
$$
\def\map#1{\smash{\mathop{\hbox to
125pt{\rightarrowfill}}\limits^{#1}}} 
\def\dmap#1#2{\smash{\mathop{\hbox to
30pt{\rightarrowfill}}\limits_{#2}\limits^{#1}}} 
\def\dpam#1#2{\smash{\mathop{\hbox to
30pt{\leftarrowfill}}\limits_{#2}\limits^{#1}}} 
\matrix {\Cal H}^j_x(i_*\Omega_{M\setminus K}^\x[2n+1])&\map{res}&
 {\Cal H}_x^j(\Dt_K)\\ 
\Bigm\|&& \Bigm\| \\ 
 H^{2n+1+j}(B_x\setminus K)\quad\dmap{\delta}{}
& H^{2n+2+j}(B_x,B_x\setminus K)\quad\dpam{[B_x]\cap}{\simeq}& 
H^{BM}_{-j}(B_x\cap K)\,.\endmatrix $$

Let ${\Cal O}^{(n+1)}_1$ be the sheaf of meromorphic forms of the type 
$(n+1,0)$
with poles of order 1 on $K$. Its sections are described locally by the 
formula
${g \over s}dz_1\wedge \dots \wedge dz_{n+1}$, where $s$ defines $K$ and $g\in
{\Cal O}_M$. In the previous chapter we have constructed a conelike metric, 
for which (under the assumption $\kappa>1$) Leray residue forms have
$p$--integrable norms:
$$Res\left({\Cal O}^{(n+1)}_1\right)\subset {\Cal L}^n_{(p)}\subset
\Omega^n_{K\setminus\Sigma}\,.$$ 
The sheaf of $L_p$--cohomology is isomorphic to an appropriate 
intersection homo\-lo\-gy sheaf \cite{We1}, so in this way we have constructed

morphisms
of sheaves: 
$${\Cal O}^{(n+1)}_1 \longrightarrow {\Cal L}^\x_{(2)} \simeq 
IC^\x_{\underline m}
\qquad \text{ and} \qquad
{\Cal O}^{(n+1)}_1 \longrightarrow {\Cal L}^\x_{(p)} \simeq IC^\x_{\underline 
0}$$
for $p>2+{2\over {n-1}}$.
A question arises: are these morphisms independent on the metric?
To be precise, let us consider the sequence of the canonical morphisms and the
obstruction sheaves \cite{GM, \S 5.5}:
$$\matrix\Ct_K[n] & \longrightarrow & IC^\x_{\underline 0}& \longrightarrow & 
IC^\x_{\underline m}
& \longrightarrow & IC^\x_{\underline t} & \longrightarrow &\Dt_K^\x\\ 
_{+1} \nwarrow & & \swarrow\,\,_{+1} \nwarrow & & \swarrow \,\,
_{+1} \nwarrow & & \swarrow\,\,_{+1} \nwarrow & & \swarrow\\ 
&S_1&&S_2&&S_3&&S_4&\\ \endmatrix \,.$$
The triangles in the diagram are distinguished in the derived category.
We regard these sheaves as sheaves on $M$ supported by $K$. The cohomology of
the links is nonzero only in dimensions 0, $n-1$, $n$ and $2n-1$, so for
arbitrary perversity the sheaf $IC^\x_{\underline p}$ is isomorphic to:
\settabs\+ ****&******&***&******************\cr
\+ {1)}& $IC_{\underline 0}$ &if& $p(2n) <n-1$,\cr
\+ {2)}& $IC_{\underline m}$ &if& $p(2n) =n-1$,\cr
\+  {3)}& $IC_{\underline t}$ &if& $p(2n) >n-1$.\cr
\noindent The obstruction sheaves $S_2$ and $S_3$ are supported by the 
singular points
and 
$$\aligned {\Cal H}_x^{-n-1}(S_2)=IH^{\underline
m}_{n+1}(cL_x)=H_n(L_x) &\quad\text{and}\quad {\Cal 
H}_x^i(S_2)=0 \quad \text{ for}\quad i\not= -(n+1)\\
 {\Cal H}_x^{-n}(S_3)=IH^{\underline 
t}_{n}(cL_x)=H_{n-1}(L_x) &\quad\text{and} \quad {\Cal 
H}_x^i(S_3)=0 \quad \text{ for}\quad i\not= -n\,.\endaligned$$ 
The obstruction sheaves $S_1$ and $S_4$ are also supported by $\Sigma_K$ and
concentrated in one dimension:
$$\aligned {\Cal H}_x^{2n}(S_0)=\tilde H^0(L_x)&\qquad\text{and}\qquad {\Cal 
H}_x^i(S_0)=0 \qquad \text{ for}\qquad i\not=-2n\\ 
 {\Cal H}_x^{0}(S_4)=\tilde H_0(L_x)&\qquad\text{and}\qquad {\Cal 
H}_x^i(S_4)=0 \qquad \text{ for}\qquad i\not=0\endaligned$$ 
Applying the functor $RHom({\Cal O}^{(n)}_1[n],-)$ to the diagram
above we get distinguished triangles
and long exact sequences. Replacing $R^0Hom$ by $Hom_\sd$ --- homomorphisms in
the derived category we obtain:
$$Hom_\sd ({\Cal O}^{(n)}_1[n],\Ct_K[2n]) @>\simeq>>
Hom_\sd ({\Cal O}^{(n)}_1[n],IC^\x_{\underline 0})$$
$$\bigoplus_{x\in\Sigma} Hom({\Cal O}^{(n)}_{1,x},H_{n}(L_x)) 
\longrightarrow
 Hom_\sd ({\Cal O}^{(n)}_1[n],IC^\x_{\underline 0}) @>\roman{epi}>>
Hom_\sd ({\Cal O}^{(n)}_1[n],IC^\x_{\underline m}) $$ 
$$ Hom_\sd ({\Cal O}^{(n)}_1[n],IC^\x_{\underline m}) 
@>\roman{mono}>>
Hom_\sd ({\Cal O}^{(n)}_1[n],IC^\x_{\underline t}) \longrightarrow{}
\bigoplus_{x\in\Sigma} Hom({\Cal O}^{(n)}_{1,x},H_{n-1}(L_x))
$$ 
$$Hom_\sd ({\Cal O}^{(n)}_1[n],IC^\x_{\underline t}) @>\simeq>> Hom_\sd ({
\Cal O}^{(n)}_1[n],\Dt^\x_K)\,.$$ 
In this way we see that:

\proclaim{Proposition 5.1} The lift of the residue morphism to ${\Cal
O}^{(n)}_1[n]\longrightarrow IC^\x_{\underline m}$ is unique. 
If such a lift exists then there exists a lift to $\Ct_K[2n]$, which is not
unique in general.
\endproclaim

The Proposition 5.1. is not a surprise since on the cohomology level we have
$$IH_n^{\underline m}(K) =
im\left(PD:\,H^n(K)\longrightarrow H_n(K)\right)$$ for $n>1$. We do not know 
if the lift to 
$IC^\x_{\underline 0}$ essentially depends on the choice of a metric. We
remind that a metric depends on the choice of coordinates in which the
singularity is quasihomogeneous. The metric was determined by the weights.

\example{Example 5.2} Consider the polynomial $s(x,y)=xy+y^{100}+z^2+t^2$ it 
is quasihomogeneous 
with weights $\frac {99}{100}$, $\frac 1{100}$, $\frac 12$ and $\frac 12$.
This is a
Morse singularity (i.e. of type $A_1$), and one can change coordinates so that
$s(x',y')={x'}^2+{y'}^2+z^2+t^2$. Then 
all weights are $\frac 12$.
\endexample

\head 6. Appendix: the $P_8$ singularity \endhead

Consider a singularity of type $P_8$:
$$s(z_1,z_2,z_3)=z_1^3+z_2^3+z_3^3$$
and let $$\omega={1\over s}dz_1\wedge dz_2\wedge dz_3\,.$$
Then $$\r= {1\over 3 z_1^2} dz_2\wedge dz_3$$
for $z_1\neq 0$. We want to check if $\r_{|L}=0\in H^2(L)$, where $L=S^5\cap
K$. The radius of the sphere does not matter as $s$ is 
homogeneous. To
calculate the cohomology of $L$ we apply the Gysin exact sequence of the 
fibration
$$S^1\hookrightarrow L@>p>> L/S^1 \subset \Pt^2\,.$$
The projectivization $L/S^1 $ of $L$ is a cubic curve in the projective plane
$\Pt^2$, so it is a topological 2--dimensional torus. We obtain the sequence:
$$@>>> H^0(L/S^1 ) @>\cup e>> H^2(L/S^1 ) @>p^*>> H^2(L) @>\int_p>> H^1(L/S^1)
@>\cup e>> H^3(L/S^1 )=0\,,$$
where the morphism $\int_p$ is the integration along the fibers of the 
projection
$p$. The bundle $L@>p>> L/S^1 $ is the restriction of the tautological bundle
$S^5 @>>> \Pt^2$. Thus the Euler class of $p$ is the restriction of the
generator of $H^2(\Pt^2)$. Hence the evaluation of the Euler class 
$\langle e,[L/S^1 ]\rangle=deg\,s=3$. Thus rationally
$\int_p$  in  the  Gysin  sequence  is  an  isomorphism  and  the 
necessary and
sufficient condition to lift is vanishing of $\left[\int_p
Res\,\omega_{|L}\right]\in H^1(L/S^1)$. We will show that this element does 
not vanish. Let 
$$U_1=\left\{[z_1:z_2:z_3]\in \Pt^2:\,z_1\neq 0\right\}
=\left\{[1:y_2:y_3]\in \Pt^2:\,y_2,\,y_3\in\Ct\right\}\simeq\Ct^2\,.$$
The tautological bundle $\tilde p:\Ct^3\setminus\{0\}@>>>\Pt^2$ restricted to
$U_1$ is trivial:
$$\aligned \tilde p^{-1}(U_1) &\simeq \Ct^*\times \Ct^2\\
(z_1,z_2,z_3)&\rightarrowtail \left(z_1,\left({z_2\over z_1},{z_3\over
z_1}\right)\right)\\ 
(y_1,y_1y_2,y_1y_3)&\leftarrowtail (y_1,y_2,y_3)\endaligned$$
We write r in $y$--coordinates:
$$r={1\over 3y_1^2}(y_1dy_2+y_2dy_1)\wedge (y_1dy_3+y_3dy_1)=
{dy_1\over 3y_1}\left(y_2dy_3-y_3dy_2\right)+{1\over 3} (dy_2\wedge dy_3)\,.$$
We integrate it over each fiber
$$p^{-1}([1,y_2,y_3])=\left\{(y_1,y_1y_2,y_1y_3):
|y_1|^2(1+|y_2|^2+|y_3|^3)=1\right\} \,;$$
$$\zeta=\int_p\r=\int_p \left[{dy_1\over 
3y_1}\left(y_2dy_3-y_3dy_2\right)+{1\over
3} (dy_2\wedge dy_3) \right]=
{\frac23 \pi i}(y_2dy_3-y_3dy_2)\,.$$
The form $\zeta$ does not vanish on $L/S^1$ since
$$\zeta\wedge d(1+y_2^3+y_3^3)={\frac23 \pi i}(y_2dy_3-y_3dy_2)\wedge
3(y_2^2dy_2+y_3^3dy_3) = -{\frac23 \pi i}(y_2^3+y_3^3)dy_2\wedge dy_3\,,$$
and it is equal ${\frac23 \pi i}\,dy_2\wedge dy_3$ on $L/S^1$.
It is a harmonic form, so its class in coho\-mo\-lo\-gy is nontrivial.
The conclusion is that $Res\,\omega\notin 
ker\left(k^*:H^n(K^\circ)\longrightarrow H^n(\partial K^\circ)\right)
$, so it has no lift to intersection homology.

\remark{Remark} The method of the example can be used to show
that
$$\left[Res\left(\frac 1sdz_1\wedge\dots\wedge dz_{n+1}\right)\right]=0\in
H^{n}(K\setminus\{0\})$$ for quasihomogeneous polynomial with
$\kappa\not=1$.  For an arbitrary $\omega=\frac gs
dz_1\wedge\dots\wedge dz_{n+1}$ one can define an obstruction in
$H^{n-1}(L/S^1)$ vanishing if and only if the $1-\kappa$ weighted homogeneous
part of $g$ vanishes. This obstruction vanishes if and only if the
residue class lifts to intersection homology \cite{We2}.
The method of obstruction does not give a lift to
cohomology. It only shows that there exists one.\endremark

\enddocument
\end